\documentclass[a4paper,english,11pt,dvips]{article}
\usepackage{citehack}
\usepackage{amsmath}
\usepackage{amsfonts}
\usepackage{amssymb}
\usepackage{longtable}
\usepackage{amsthm}
\usepackage{euscript}
\usepackage{array}

\newcommand{\REM}[1]{\relax}
\makeatletter%

\newcommand{\begthrm}[2]{\begin{trivlist}\it\item[\hspace{\labelsep}{\bf #1~#2.}]}
\newcommand{\enthrm}{\end{trivlist}}

\renewcommand{\@makefntext}[1]{{\normalsize\parindent=1em\noindent\hbox to 1.8em{\hss
$^{\@thefnmark}$}#1}}

\renewcommand{\p@equation}{\relax}

\renewcommand{\l@subsection}{\@dottedtocline{3}{1.5em}{3.3em}}

\expandafter\def\csname @tocrmarg\endcsname{4em}

\thm@notefont{\rm}

\newtheorem{theorem}{Theorem}

\newtheorem{definit}{Definition}

\newtheorem{ex}{Example}

\newtheorem{remark}{Remark}

\def\ZR{\rightthreetimes}

\def\spa{\mathop\text{{\rm span}}\nolimits}

\def\Hom{\mathop\text{\rm Hom}\nolimits}

\def\rk{\mathop\text{\rm rk}\nolimits}
\def\pr{\mathop\text{\rm pr}\nolimits}

\def\Aut{\mathop\text{\rm Aut}\nolimits}
\def\End{\mathop\text{\rm End}\nolimits}
\def\real{\mathbb{R}}

\def\k{\mathfrak{k}}
\def\g{\mathfrak{g}}
\def\h{\mathfrak{h}}
\def\so{\mathfrak{so}}

\def\u{\mathfrak{u}}

\def\sp{\mathfrak{sp}}

\def\Simil{\mathop\text{\rm Sim}\nolimits}

\def\A{\mathcal {A}}
\def\B{\mathcal {B}}

\def\N{\mathcal {N}}

\def\id{\mathop\text{\rm id}\nolimits}

\def\im{\mathop\text{\rm Im}\nolimits}
\def\re{\mathop\text{\rm Re}\nolimits}
\def\lie{\mathop{\mathcal{LA}}\nolimits}
\def\Mat{\mathop\text{\rm Mat}\nolimits}
\def\Op{\mathop\text{\rm Op}\nolimits}
\def\Sp{\mathop\text{\rm Sp} \nolimits}
\def\SO{\mathop\text{\rm SO}\nolimits}
\def\U{\mathop\text{\rm U}\nolimits}

\title{Weakly irreducible subgroups of $\Sp(1,n+1)$}

\author{Natalia I. Bezvitnaya}

\begin{document}
\maketitle\vskip-50ex {\renewcommand{\abstractname}{Abstract}\begin{abstract} Connected weakly irreducible not irreducible
subgroups of $\Sp(1,n+1)\subset\SO(4,4n+4)$ that satisfy a certain additional condition are classified. This will be used
to classify connected holonomy groups of pseudo-hyper-K\"ahlerian manifolds of index 4.
\end{abstract}


{\bf Mathematical subject codes:}  53C29, 53C50


\section{Introduction} The classification of connected holonomy groups of Riemannian manifolds is well known \cite{Ber,Be,Bryant,Jo}.
A classification of holonomy groups of pseudo-Riemannian manifolds is an actual problem of differential geometry.
 Very recently were obtained classifications of connected holonomy groups of Lorentzian manifolds
 \cite{B-I,Le,Gal5} and of pseudo-K\"ahlerian manifolds of index 2 \cite{GalKahler}. These groups are contained in
 $\SO(1,n+1)$ and $\U(1,n+1)\subset\SO(2,2n+2)$, respectively. As the next step, we study connected holonomy groups
 contained in $\Sp(1,n+1)\subset\SO(4,4n+4)$, i.e. holonomy groups of pseudo-hyper-K\"ahlerian manifolds of index 4.
By the Wu theorem \cite{Wu} and the results of Berger for connected  irreducible holonomy groups of pseudo-Riemannian
manifolds \cite{Ber}, it is enough to consider only weakly irreducible not irreducible holonomy groups (each such group
does not preserve any proper non-degenerate vector subspace of the tangent space, but preserves a degenerate subspace).

In the present paper we classify connected weakly irreducible not irreducible subgroups of $\Sp(1,n+1)\subset\SO(4,4n+4)$
($n\geq 1$) that satisfy a natural condition. The case $n=0$ will be considered separately. We generalize the method of
\cite{Gal1,GalKahler}. Let $G\subset\Sp(1,n+1)$ be a weakly irreducible not irreducible subgroup and $\g\subset\sp(1,n+1)$
the corresponding subalgebra. The results of \cite{GalKahler} allow us to expect that if $\g$ is the holonomy algebra,
then $\g$ containes a certain $3$-dimensional ideal $\B$. We will prove this in another paper. Consider the action of $G$
on the space $\mathbb{H}^{1,n+1}$, then $G$ acts on the boundary of the quaternionic hyperbolic space, which is
diffeomorphic to the $4n+3$-dimensional sphere $S^{4n+3}$ and $G$ preserves a point of this space. We define a map
$s_1:S^{4n+3}\backslash\{point\}\to\mathbb{H}^n$ similar to the usual stereographic projection. Then any $f\in G$ defines
the map $F(f)=s_1\circ f\circ s_2:\mathbb{H}^n\to\mathbb{H}^n$, where $s_2:\mathbb{H}^n\to S^{4n+3}\backslash\{point\}$ is
the inverse of the usual stereographic projection restricted to
$\mathbb{H}^n\subset\mathbb{H}^n\oplus\mathbb{R}^3=\mathbb{R}^{4n+3}$. We get that $F(G)$ is contained in the group
$\Simil \mathbb{H}^n$ of similarity transformations of $\mathbb{H}^n$. We show that $F(G)$ preserves an affine subspace
$L\subset\mathbb{R}^{4n}=\mathbb{H}^{n}$ such that the minimal affine subspace of $\mathbb{H}^{n}$ containing $L$ is
$\mathbb{H}^{n}$. Moreover, $F(G)$ does not preserve any proper affine subspace of $L$. Then $F(G)$ acts transitively on
$L$ \cite{Al2}. We describe subspaces $L$ with such property and using results of \cite{GalKahler} we find all connected
Lie subgroups $K \subset \Simil \mathbb{H}^{n}$ preserving $L$ and acting transitively on $L$. Note that the kernel of the
Lie algebra homomorphism $dF:\g\to\lie(\Simil \mathbb{H}^{n})$ coincides with the ideal $\B$. Consequently,
$\g=(dF)^{-1}(\k)$, where $\k\subset\lie(\Simil \mathbb{H}^{n})$ is the Lie algebra of one of the obtained Lie subgroups
$K \subset \Simil \mathbb{H}^{n}$.

Note that we classify weakly irreducible not irreducible subgroups of $\Sp(1,n+1)$ up to conjugacy in $\SO(4,4n+4)$. It is
also possible to classify these subgroups up to conjugacy in $\Sp(1,n+1)$, see Remark \ref{remark1}.

{\it Acknowledgement.} I am grateful to Jan Slov\'ak for support and help. The author has been supported by the grant GACR
201/05/H005.

\section{Preliminaries}

First we summarize some facts about quaternionic vector spaces. Let $\mathbb{H}^{m}$ be an m-dimensional quaternionic
vector space and $e_{1},...,e_{m}$ a basis of $\mathbb{H}^{m}$. We identify an element $X\in\mathbb{H}^{m}$ with the
column $(X_{t})$ of the left coordinates of $X$ with respect to this basis, $X=\sum_{t=1}^{m}X_{t}e_{t}$.

Let $f:\mathbb{H}^{m}\to\mathbb{H}^{m}$ be an $\mathbb{H}$-linear map. Define the matrix $\Mat_{f}$ of $f$ by the relation
$fe_{l}=\sum_{t=1}^{m}(\Mat_{f})_{tl}e_{t}$. Now if $X\in\mathbb{H}^{m}$, then $fX=(X^{t}\Mat_{f}^{t})^{t}$ and because of
the non-commutativity of the quaternions this is not the same as $\Mat_{f}X$. Conversely, to an $m\times m$ matrix $A$ of
the quaternions we put in correspondence the linear map $\Op{A}:\mathbb{H}^{m}\to\mathbb{H}^{m}$ such that $\Op{A}\cdot
X=(X^{t}A^{t})^{t}$. If $f,g:\mathbb{H}^{m}\to\mathbb{H}^{m}$ are two $\mathbb{H}$-linear maps, then
$\Mat_{fg}=(\Mat_{g}^{t}\Mat_{f}^{t})^{t}$. Note that the multiplications by the imaginary quaternions are not
$\mathbb{H}$-linear maps. Also, for $a,b\in\mathbb{H}$ holds $\overline{ab}=\bar{b}\bar{a}$. Consequently, for two square
quaternionic matrices we have $(\overline{AB})^{t}=\bar{B}^t\bar{A}^t$.

A pseudo-quaternionic-Hermitian metric $g$ on $\mathbb{H}^{m}$ is a non-degenerate $\mathbb{R}$-bilinear map
$g:\mathbb{H}^{m}\times\mathbb{H}^{m}\to\mathbb{H}$ such that $g(aX,Y)=ag(X,Y)$ and $\overline{g(Y,X)}=g(X,Y)$, where
$a\in\mathbb{H}$, $X,Y\in\mathbb{H}^{m}$. Hence, $g(X,aY)=g(X,Y)\bar{a}$. There exists a basis $e_{1},...,e_{m}$ of
$\mathbb{H}^{m}$ and integers $(r,s)$ with $r+s=m$ such that $g(e_{t},e_{l})=0$ if $t \neq l$, $g(e_{t},e_{t})=-1$ if $1
\leq t \leq p$ and $g(e_{t},e_{t})=1$ if $p+1 \leq t \leq m$. The pair $(r,s)$ is called  the signature of $g$. In this
situation we denote $\mathbb{H}^{m}$ by $\mathbb{H}^{r,s}$. The realification of $\mathbb{H}^{m}$ gives us the vector
space $\mathbb{R}^{4m}$ with the quaternionic structure $(i,j,k)$. Conversely, a quaternionic structure on
$\mathbb{R}^{4m}$, i.e. a triple $(I,J,K)$ of endomorphisms of $\mathbb{R}^{4m}$ such that $I^{2}=J^{2}=K^{2}=-\id$ and
$K=IJ=-JI$, allows us to consider $\mathbb{R}^{4m}$ as $\mathbb{H}^{m}$. A pseudo-quaternionic-Hermitian metric $g$ on
$\mathbb{H}^{m}$ of signature $(r,s)$ defines on $\mathbb{R}^{4m}$ the $i,j,k$-invariant pseudo-Euclidian metric $\eta$ of
signature $(4r,4s)$, $\eta(X,Y)=\re g(X,Y)$, $X,Y\in\mathbb{R}^{4m}$. Conversely, a $I,J,K$-invariant pseudo-Euclidian
metric on $\mathbb{R}^{4m}$ defines a pseudo-quaternionic-Hermitian metric $g$ on $\mathbb{H}^{m}$,
$$g(X,Y)=\eta(X,Y)+i\eta(X,IY)+j\eta(X,JY)+k\eta(X,KY).$$ The Lie group $\Sp(r,s)$ and its Lie algebra $\sp(r,s)$ are
defined as follows \begin{align*}\Sp(r,s)=&\{f\in\Aut(\mathbb{H}^{r,s})|\,\, g(fX,fY)=g(X,Y) \text{ for all }
X,Y\in\mathbb{H}^{r,s}\},\\ \sp(r,s)=&\{f\in\End(\mathbb{H}^{r,s})|\, g(fX,Y)+g(X,fY)=0 \text{ for all }
X,Y\in\mathbb{H}^{r,s}\}.\end{align*}

\section{The Main Theorem}\label{MainTh}

\begin{definit} A subgroup $G\subset \SO(r,s)$ (or a subalgebra $\g\subset\so(r,s)$) is called weakly irreducible if it does
not preserve any non-degenerate proper vector subspace of $\real^{r,s}$. \end{definit}

Let $\real^{4,4n+4}$ be a $(4n+8)$-dimensional real vector space endowed with a quaternionic structure $I,J,K\in
\End(\real^{4,4n+4})$ and  an $I, J, K$-invariant metric $\eta$ of signature $(4,4n+4)$. We identify this space with the
$(n+2)$-dimensional quaternionic  space $\mathbb{H}^{1,n+1}$ endowed  with the pseudo-quaternionic-Hermitian metric $g$ of
signature $(1,n+1)$ as above.

Obviously, if a Lie subgroup $G\subset \Sp(1,n+1)$ acts weakly irreducibly not irreducibly on $\real^{4,4n+4}$, then $G$
acts weakly irreducibly not irreducibly on $\mathbb{H}^{1,n+1}$. The converse is not true, see Example \ref{example2}
below. If $G$ acts weakly irreducibly not irreducibly on $\mathbb{H}^{1,n+1}$, then $G$ preserves a proper degenerate
subspace $W\subset\mathbb{H}^{1,n+1}$. Consequently, $G$ preserves the intersection $W\cap
W^{\perp}\subset\mathbb{H}^{1,n+1}$, which is an isotropic quaternionic line.

Fix a Wit basis $p,e_1,...,e_n,q$ of $\mathbb{H}^{1,n+1}$, i.e. the Gram matrix of the metric $g$ with respect to this
basis has the form $\left (\begin{array}{ccc} 0 & 0 & 1\\ 0 & E_n & 0 \\ 1 & 0 & 0 \\
\end{array}\right),$ where $E_n$ is the $n$-dimensional identity matrix.
Denote by $\Sp(1,n+1)_{\mathbb{H}{p}}$ the Lie subgroup of $\Sp(1,n+1)$ acting on $\mathbb{H}^{1,n+1}$ and preserving the
quaternionic isotropic line $\mathbb{H}{p}$. Note that any weakly irreducible and not irreducible subgroup of $\Sp
(1,n+1)$ is conjugated to a weakly irreducible subgroup of $\Sp(1,n+1)_{\mathbb{H}{p}}$. The Lie subalgebra
$\sp(1,n+1)_{\mathbb{H}{p}}\subset \sp(1,n+1)$ corresponding to the Lie subgroup\\ $\Sp(1,n+1)_{\mathbb{H}{p}}\subset \Sp
(1,n+1)$ has the following  form $$\sp(1,n+1)_{\mathbb{H}{p}}=\left\{\left. \Op\left(\begin{array}{cccc} \bar{a} &-\bar
X^t & b\\ 0 & \Mat_{h} &X
\\ 0 & 0 & -a \\
\end{array}\right)\right|\begin{matrix} a\in \mathbb{H},\,\,  X\in \mathbb{H}^{n},\\
  h\in \sp(n),\,\,b\in\im\mathbb{H} \end{matrix}
\right\}.$$ Let $(a,A,X,b)$ denote the above element of $\sp(1,n+1)_{\mathbb{H}{p}}$. Define the following vector
subspaces of $\sp(1,n+1)_{\mathbb{H}{p}}:$
\begin{align*}\A_{1}=&\{(a,0,0,0)| a\in\mathbb{R}\},& \A_{2}=&\{(a,0,0,0)| a\in\im\mathbb{H}\},\\ \N=&\{(0,0,X,0)|
X\in\mathbb{H}^{n}\},& \B=&\{(0,0,0,b)| b\in\im\mathbb{H}\}.\end{align*} Obviously, $\sp(n)$ is a subalgebra of
$\sp(1,n+1)_{\mathbb{H}{p}}$ with the inclusion $$h\in\sp(n)\mapsto
\Op\left(\begin{matrix}0&0&0\\0&\Mat_h&0\\0&0&0\end{matrix}\right)\in\sp(1,n+1)_{\mathbb{H}{p}}.$$ We obtain that $\A_{1}$
is a one-dimensional commutative subalgebra that commutes with $\A_2$ and $\sp(n)$, $\A_{2}$ is a subalgebra isomorphic to
$\sp(1)$ and commuting with $\sp(n)$, $\B$ is a commutative ideal, which commutes with $\sp(n)$ and $\N$. Also,
\begin{align*}[(a,0,0,0),(0,0,X,b)]=&(0,0,aX,2\im ab),\quad [(0,0,X,0),(0,0,Y,0)]=(0,0,0,2 \im g(X,Y)),\\
[(0,A,0,0),(0,0,X,0)]=&(0,0,(X^{t}A^{t})^{t},0),\end{align*}
 where $a\in \mathbb{H}$, $X,Y\in \mathbb{H}^{n}$, $A=\Mat_{h}$, $ h\in
\sp(n)$, $b\in \im\mathbb{H}$. Thus we have the decomposition
$$\sp(1,n+1)_{\mathbb{H}{p}}=(\A_{1}\oplus\A_{2}\oplus\sp(n))\ltimes(\N+\B)\simeq(\mathbb{R}\oplus\sp(1)\oplus\sp(n))\ltimes(\mathbb{H}^{n}+\mathbb{R}^{3}).$$
Now consider two examples.

\begin{ex}\label{example1} The subalgebra $\g=\{(0,0,X,b)|\,\, X\in\mathbb{R}^{n},\, b\in\im\mathbb{H}\}\subset\sp(1,n+1)_{\mathbb{H}{p}}$ acts weakly irreducibly on $\mathbb{R}^{4,4n+4}$.
\end{ex}
{\it Proof.} Assume the converse. Let $\g$ preserve a non-degenerate proper vector subspace $L\subset
\mathbb{R}^{4,4n+4}$. Suppose  the projection of $L$ to $\mathbb{H}{q}\subset\mathbb{H}^{1,n+1}=\mathbb{R}^{4,4n+4}$ is
non-zero, then there is a vector $v\in L$ such that $v=v_0p+v_1+v_2q$, where $v_0, v_2\in \mathbb{H}$, $v_2\neq0$ and
$v_1\in\mathbb{H}^{n}$. Consider  elements $\xi_1=(0,0,X,0)\in \g$ with $g(X,X)=1$ and $\xi_2=(0,0,0,b)\in \g$. Then,
$\xi_1(\xi_1v)=-v_2p\in L$ and $\xi_2v=v_2bp\in L$. Since $v_2\neq0$, we have $\mathbb{H}{p}\subset L$. It follows that
$L^{\perp_{\eta}}\subset \mathbb{H}{p}\oplus \mathbb{H}^{n}$ and $L^{\perp_{\eta}}$ is a $\g$-invariant non-degenerate
proper subspace. Now we can assume that $\g$ preserves a non-trivial non-degenerate vector subspace $L\subset
\mathbb{H}{p}\oplus \mathbb{H}^{n}$. Let $v=v_0p+v_1\in L$, $v\neq 0$. If $v_1=0$, then $L$ is degenerate. If $v_1\neq 0$,
then there is $X\in \mathbb{R}^{n}$ with $g(v_1,X)\neq 0$. We get $(0,0,X,0)v=-g(v_1,X)p\in L$. Hence $L$ is degenerate.
Thus we have a contradiction. $\Box$

\begin{ex}\label{example2} The subalgebra $\g=\{(0,0,X,0)|\,\, X\in\mathbb{R}^{n}\}\subset\sp(1,n+1)_{\mathbb{H}{p}}$
acts weakly irreducibly on $\mathbb{H}^{1,n+1}$ and not weakly irreducibly on $\mathbb{R}^{4,4n+4}$.
\end{ex}
{\it Proof.} The proof of the first statement is similar to the proof of Example \ref{example1}.
Clearly, the subalgebra $\g$ preserves the non-degenerate vector subspace $\spa_{\mathbb{R}}\{p,e_1,...,e_n,q\}\subset
\mathbb{R}^{4,4n+4}$. $\Box$

The classification of the holonomy algebras contained in $\u(1,n+1)$ \cite{GalKahler} gives us the following hypothesis:
{\it If $n\geq1$ and $\g\subset\sp(1,n+1)_{\mathbb{H}{p}}$ is a holonomy algebra, then $\g$ containes the ideal $\B$.} We
will prove this hypothesis in an other paper.

In the following theorem we denote the real vector subspace  $L\subset\mathbb{R}^{4n}=\mathbb{H}^{n}$ of the form
$$L={\spa}_{\mathbb{H}}\{e_{1},...,e_{m}\}\oplus {\spa}_{\mathbb{R}\oplus i\mathbb{R}}\{e_{m+1},...,e_{m+k}\}\oplus {\rm
span}_{\mathbb{R}}\{e_{m+k+1},...,e_{n}\}$$  by $\mathbb{H}^{m}\oplus \mathbb{C}^{k}\oplus \mathbb{R}^{n-m-k}$. Let
$\u(k)$ be the subalgebra of $\sp({\rm span}_{\mathbb{H}}\{e_{m+1},...,e_{m+k}\})$ that consists of the elements
$\Op\left(\begin{matrix}A&0\\0&A\end{matrix}\right)$, where $A\in\u({\spa}_{\mathbb{R}\oplus
i\mathbb{R}}\{e_{m+1},...,e_{m+k}\})$ and we use the decomposition ${\spa}_{\mathbb{H}}\{e_{m+1},...,e_{m+k}\}={\rm
span}_{\mathbb{R}\oplus i\mathbb{R}}\{e_{m+1},...,e_{m+k}\}+j{\spa}_{\mathbb{R}\oplus
i\mathbb{R}}\{e_{m+1},...,e_{m+k}\}.$ Similarly, let $\so(n-m-k)$ be the subalgebra of $\sp({\rm
span}_{\mathbb{H}}\{e_{m+k+1},...,e_{n}\})$ that consists of the elements
$\Op\left(\begin{matrix}A&0&0&0\\0&A&0&0\\0&0&A&0\\0&0&0&A\end{matrix}\right)$, where $A\in\so({\rm
span}_{\mathbb{R}}\{e_{m+k+1},...,e_{n}\})$ and we use the decomposition $\mathbb{H}^{n-m-k}=\mathbb{R}^{n-m-k}\oplus
i\mathbb{R}^{n-m-k}\oplus j\mathbb{R}^{n-m-k}\oplus k\mathbb{R}^{n-m-k}.$ For a Lie algebra $\h$ we denote by $\h'$ the
commutant $[\h,\h]$ of $\h$.

\begin{theorem}  Let $n\geq1$. Any weakly irreducible subalgebra of $\sp(1,n+1)_{\mathbb{H}{p}}$ that
contains the ideal $\B$ is conjugated by an element of $\SO(4,4n+4)$ to one of the following subalgebras:

\begin{description}

\item{{\bf Type I.}} $\g=\{(a_1+a_2,A,X,b)|\,\, a_1\in\real,\,\, a_2\in\h_0,\,\, A\in\h,\,\,X\in\mathbb{H}^{n},\,\,b\in\im\mathbb{H}\},$
where $\h_0\subset\sp(1)$ is a subalgebra of dimension $2$ or $3$, $\h\subset\sp(n)$ is a subalgebra.

\item{{\bf Type II.}} $\g=\{(a_1+ta_2+\phi(A),A,X,b)|\,a_1,t\in\real,\,\, A\in\h,\,\, X\in\mathbb{H}^{n},\,\, b\in
{\im}\mathbb{H}\},$ where $a_2\in\sp(1)$, $\h\subset\sp(n)$ is a subalgebra, $\phi:\h\to\sp(1)$ is a homomorphism.\\ If
$a_2\neq0$, then $\rk  \phi\leq1$ and $[\im\phi,a_2]\subset\mathbb{R}a_2$.

\item{{\bf Type III.}} $\g=\{(\varphi(a_2,A)+a_2,A,X,b)|\, a_2\in\h_0,\,\, A\in\h,\,\, X\in\mathbb{H}^{n},\,\, b\in{\im}\mathbb{H}\},$ where $\h_0\subset\sp(1)$ is a subalgebra of dimension $2$ or $3$, $\h\subset\sp(n)$ is a subalgebra,
$\varphi\in\Hom(\h_0\oplus\h,\mathbb{R})$, $\varphi|_{\h_0'\oplus\h'}=0$.  In particular, if $\dim\h_0=3$, i.e.
$\h_0=\sp(1)$, then $\varphi|_{\h_0}=0$.

\item{{\bf Type IV.}} $\g=\{(\varphi(t,A)+ta_2+\phi(A),A,X,b)|\, t\in\mathbb{R},\,\, A\in\h,\,\, X\in\mathbb{H}^{n},\,\, b\in\im \mathbb{H}\},$ where $a_2\in\sp(1)$, $\h\subset\sp(n)$ is a subalgebra,
$\varphi\in \Hom(\mathbb{R}\oplus\h,\mathbb{R})$, $\varphi|_{\h'}=0$, $\phi:\h\to\sp(1)$ is a homomorphism. If $a_2\neq0$,
then $\rk \phi\leq1$ and $[\im\phi,a_2]\subset\mathbb{R}a_2$. If $a_2\neq0$ and $\phi\neq0$, then
$\varphi|_{\mathbb{R}}=0$.

\item{{\bf Type V.}} $\g=\{(a_1+a_2i,A,X,b)|\, a_1,a_2\in\mathbb{R},\,\, A\in\h,\,\, X\in\mathbb{H}^{m}\oplus\mathbb{C}^{n-m},
\,\, b\in
\im\mathbb{H}\},$ where $0\leq m<n$, $\h\subset\sp(m)\oplus\u(n-m)$ is a subalgebra.

\item{{\bf Type VI.}} $\g=\{(a_1+\phi(A)i,A,X,b)|\, a_1\in\mathbb{R},\,\, A\in\h,\,\, X\in\mathbb{H}^{m}\oplus\mathbb{C}^{k}\oplus\mathbb{R}^{n-m-k},\,\, b\in\im\mathbb{H}\},$ where $0\leq m<n$, $0\leq k\leq n-m$, $\h\subset\sp(m)\oplus\u(k)\oplus\so(n-m-k)$ is a subalgebra,
$\phi\in\Hom(\h,\mathbb{R})$, $\phi|_{\h'}=0$. If $n-m-k\geq1$, then $\phi=0$.

\item{{\bf Type VII.}} $\g=\{(\varphi(a_2,A)+a_2i,A,X,b)|\, a_2\in\mathbb{R},\,\, A\in\h,\,\, X\in\mathbb{H}^{m}\oplus\mathbb{C}^{n-m},\,\, b\in\im\mathbb{H}\},$ where $0\leq m<n$, $\h\subset\sp(m)\oplus\u(n-m)$ is a subalgebra, $\varphi\in{\rm
Hom}(\mathbb{R}\oplus\h,\mathbb{R})$, $\varphi|_{\h'}=0$.

\item{{\bf Type VIII.}} $\g=\{(\varphi(A)+\phi(A)i,A,X,b)|\, A\in\h,\,\, X\in\mathbb{H}^{m}\oplus\mathbb{C}^{k}\oplus\mathbb{R}^{n-m-k},\,\, b\in\im\mathbb{H}\},$ where $0\leq m<n$, $0\leq k\leq n-m$, $\h\subset\sp(m)\oplus\u(k)\oplus\so(n-m-k)$ is a subalgebra,
$\varphi,\phi\in\Hom(\h,\mathbb{R})$, $\varphi|_{\h'}=\phi|_{\h'}=0$. If $n-m-k\geq1$, then $\phi=0$.

\item{{\bf Type IX.}} $\g=\{(0,A,\psi(A)+X,b)|\, A\in\h,\,\, X\in W,\,\, b\in\im\mathbb{H}\}.$  Here $0\leq m\leq n$ and $0\leq k\leq n-m$. For
$L=\mathbb{H}^{m}\oplus\mathbb{C}^{k}\oplus\mathbb{R}^{n-m-k}\subset\mathbb{R}^{4n}=\mathbb{H}^{n}$ we have an
$\eta$-orthogonal decomposition $L=W\oplus U$, $\h\subset\sp(W\cap iW\cap jW\cap kW)$ is a subalgebra and $\psi:\h\to W$
is a surjective linear map with $\psi|_{\h'}=0$.

\end{description}\end{theorem}

\section{Relation with the group of similarity transformations of~$\mathbb{H}^{n}$}

Let $\mathbb{H}^{n}$ be the $n$-dimensional quaternionic vector space endowed with a quaternionic-Hermitian metric $g$.
For elements $a_1\in\mathbb{R}_{+}$, $a_2\in\Sp(1), f\in\Sp(n)$ and $X\in\mathbb{H}^{n}$ consider the following
transformations of $\mathbb{H}^{n}$: $d(a_1):Y\mapsto a_1Y$ (real dilation),  $a_2:Y\mapsto a_2Y$ (quaternionic dilation),
$f:Y\mapsto fY$ (rotation), $t(Y):Y\mapsto Y+X$ (translation), here $Y\in\mathbb{H}^{n}$. Note that the elements
$a_2\in\Sp(1)$ act on $\mathbb{H}^{n}$ as $\mathbb{R}$-linear (but not $\mathbb{H}$-linear) isomorphism. These
transformations generate the Lie group $\Simil \mathbb{H}^{n}$ of similarity transformations of $\mathbb{H}^{n}$. We get
the decomposition $$\Simil \mathbb{H}^{n}=(\mathbb{R}_{+}\times\Sp(1)\cdot\Sp(n))\ZR\mathbb{H}^{n}.$$ The Lie group
$\Simil\mathbb{H}^{n}$ is a Lie subgroup of the connected Lie group $\Simil ^{0}\mathbb{R}^{4n}$ of similarity
transformations of $\mathbb{R}^{4n}$, $\Simil ^{0}\mathbb{R}^{4n}=(\mathbb{R}_{+}\times{\SO}(4n))\ZR\mathbb{R}^{4n}.$\\
The corresponding Lie algebra $\lie  (\Simil \mathbb{H}^{n})$ to the Lie group $\Simil \mathbb{H}^{n}$ has the following
decomposition $$\lie  (\Simil \mathbb{H}^{n})=(\mathbb{R}\oplus\sp(1)\oplus\sp(n))\ltimes\mathbb{H}^{n}.$$

Let $p,e_1,...,e_n,q$ be the basis of $\mathbb{H}^{1,n+1}$ as above. Consider also the basis $e_0,e_1,...,e_n,e_{n+1}$,
where $e_0=\frac{\sqrt{2}}{2}(p-q)$ and $e_{n+1}=\frac{\sqrt{2}}{2}(p+q)$. With respect to this basis the Gram matrix of
$g$ has the form $\left(\begin{matrix}-1&0\\0&E_{n+1}\end{matrix}\right).$

The subset of the $(n+1)$-dimensional quaternionic projective space $\mathbb{PH}^{1,n+1}$ that consists of all
quaternionic isotropic lines is called the {\it boundary} of the quaternionic hyperbolic space and is denoted by
$\partial\mathbf{H}^{n+1}_\mathbb{H}.$

Let $h_0,...,h_{n+1}$, where $h_s=x_s+iy_s+jz_s+kw_s\in\mathbb{H}$ ($0\leq s\leq n+1$) be the coordinates on
$\mathbb{H}^{1,n+1}$ with respect to the basis $e_0,...,e_{n+1}.$ Denote by $\mathbb{H}^{n}$ and $\mathbb{H}^{n+1}$ the
subspaces of $\mathbb{H}^{1,n+1}$ spanned by the vectors $e_1,...,e_n$ and $e_1,...,e_{n+1}$, respectively. Note that the
intersection $(e_0+\mathbb{H}^{n+1})\cap\{X\in\mathbb{H}^{1,n+1}|\, g(X,X)=0\}$ is given by the system of equations:
$$x_0=1, y_0=0, z_0=0, w_0=0, x_1^{2}+y_1^{2}+z_1^{2}+w_1^{2}+...+x_{n+1}^{2}+y_{n+1}^{2}+z_{n+1}^{2}+w_{n+1}^{2}=1,$$
i.e. this set is the $(4n+3)$-dimensional unite sphere $S^{4n+3}$. Moreover, each isotropic line intersects  this set at a
unique point, e.g. $ \mathbb{H}{p}$ intersects it at the point $\sqrt{2}p.$ Thus we identify the space
$\partial\mathbf{H}^{n+1}_\mathbb{H}$ with the sphere $S^{4n+3}$. Any $f\in\Sp(1,n+1)_{\mathbb{H}p}$ takes quaternionic
isotropic lines to quaternionic isotropic lines and preserves the quaternionic isotropic line $\mathbb{H}p$. Hence it acts
on $\partial\mathbf{H}^{n+1}_\mathbb{H}\setminus\{\mathbb{H}p\}=S^{4n+3}\setminus\{\sqrt{2}p\}$.

Consider the connected Lie subgroups $A_1, A_2,\Sp(n)$ and $P$ of $\Sp(1,n+1)_{\mathbb{H}p}$ corresponding to the
subalgebras $\A_1, \A_2,\sp(n)$ and $\N+\B$ of the Lie algebra $\sp(1,n+1)_{\mathbb{H}p}$. With respect to the basis
$p,e_1,...,e_n,q$ these groups have the following matrix form: \begin{align*} A_1=&\left\{\Op\left. \left
(\begin{array}{ccc} a_1 & 0 & 0\\ 0 & E_n & 0\\ 0 & 0 & a_1^{-1}
\end{array} \right)\right|\begin{array}{c} a_1\in \real_+
\end{array} \right \},\quad
A_2=\left\{\Op\left. \left (\begin{array}{ccc} e^{-a_2} & 0 & 0\\ 0 & E_n & 0\\ 0 & 0 & e^{-a_2}
\end{array} \right)\right| a_2\in \im \mathbb{H} \right \},\\
\Sp(n)=&\left\{\Op\left. \left (\begin{array}{ccc} 1 & 0 & 0\\ 0 & \Mat_f & 0\\ 0 & 0 & 1
\end{array} \right)\right| f\in\Sp(n) \right \},\\
P=&\left\{\Op\left. \left ( \begin{array}{ccc} 1 & -\bar{Y}^{t} & b-\frac{1}{2}Y^t\bar{Y}\\ 0 & E_n & Y\\ 0 & 0& 1
\end{array} \right)\right|\begin{array}{c} Y\in\mathbb{H}^{n},\\ b\in\im\mathbb{H} \end{array} \right \}.\end{align*}
 We have the decomposition $$\Sp(1,n+1)_{\mathbb{H}p}=(A_1\times A_2\times \Sp(n))\ZR
P\simeq(\mathbb{R}_{+}\times\Sp(1)\times\Sp(n))\ZR(\mathbb{H}^{n}\cdot\mathbb{R}^{3}).$$ Let
$s_1:S^{4n+3}\setminus\{\sqrt{2}p\}\to e_0+\mathbb{H}^{n}$ be the map defined as the usual stereographic projection, but
using quaternionic lines. More precisely, for $s\in S^{4n+3}\setminus\{\sqrt{2}p\}$ we define $s_1(s)$ to be the point of
the intersection of $e_0+\mathbb{H}^{n}$ with the quaternionic line passing through the points $\sqrt{2}p$ and $s$. It is
easy to see that this intersection consists of a single point. Let $s_2:e_0+\mathbb{H}^{n}\to
S^{4n+3}\setminus\{\sqrt{2}p\}$ be the restriction to $e_0+\mathbb{H}^{n}$ of the inverse to the usual stereographic
projection from $S^{4n+3}\setminus\{\sqrt{2}p\}$ to $e_0+\mathbb{H}^{n}\oplus(\im\mathbb{H})e_{n+1}$. Note that $s_1\circ
s_2=\id_{e_0+\mathbb{H}^{n}}$, but unlike in the usual case, $s_1$ is not surjective. We have $s_2\circ s_1|_{\im
s_2}=\id_{\im s_2}$. Also, let $e_0$ and $-e_0$ denote the translations $\mathbb{H}^{n}\to e_0+\mathbb{H}^{n}$ and
$e_0+\mathbb{H}^{n}\to \mathbb{H}^{n}$, respectively.

For $f\in\Sp(1,n+1)_{\mathbb{H}{p}}$ define the map $$F(f)=(-e_0)\circ s_1\circ f\circ s_2\circ e_0:\mathbb{H}^{n}\to
\mathbb{H}^{n}.$$ Now we will show that $F$ is a surjective homomorphism from the Lie group $\Sp(1,n+1)_{\mathbb{H}{p}}$
to the Lie group $\Simil \mathbb{H}^{n}$ and $\ker F=\mathbb{Z}_2\times B$, where $\mathbb{Z}_2=\{\id,-\id\}\in
\Sp(1,n+1)_{\mathbb{H}{p}}$ and $B$ is the connected Lie subgroup of $\Sp(1,n+1)_{\mathbb{H}{p}}$ corresponding to the
ideal $\B\subset\sp(1,n+1)_{\mathbb{H}{p}}$. First of all, the computations show that for $a_1\in\mathbb{R}$, $a_2\in\im
\mathbb{H}$, $f\in\Sp(n)$ and $Y\in\mathbb{H}^{n}$ it holds
$$F\left(\Op\left(\begin{matrix}a_1&0&0\\0&E_n&0\\0&0&a_1^{-1}\end{matrix}\right)\right)=d(a_1)\in\mathbb{R}_{+}\subset\Simil\mathbb{H}^{n},$$
$$F\left(\Op\left(\begin{matrix}e^{-a_2}&0&0\\0&E_n&0\\0&0&a^{-a_2}\end{matrix}\right)\right)=e^{a_2}\in\Sp(1)\subset\Simil\mathbb{H}^{n},$$
$$F\left(\Op\left(\begin{matrix}1&0&0\\0&\Mat_{f}&0\\0&0&1\end{matrix}\right)\right)=f\in\Sp(n)\subset\Simil\mathbb{H}^{n},$$
$$F\left(\Op\left(\begin{matrix}1&-\bar{Y}^{t}&b-\frac{1}{2}Y^{t}\bar{Y}\\0&E_n&Y\\0&0&1\end{matrix}\right)\right)=t\left(-\frac{\sqrt{2}}{2}Y\right)\in\mathbb{H}^{n}\subset\Simil\mathbb{H}^{n}.$$
It follows that if $f_1,f_2\in P$, then $F(f_1f_2)=F(f_1)F(f_2)$, i.e. $F|_{P}$ is a homomorphism from $P$ to
$\Simil\mathbb{H}^{n}$. It can easily be checked that any $f\in A_1\times A_2\times\Sp(n)$ considered as a map from
$S^{4n+3}\setminus\{\sqrt{2}p\}$ to itself preserves $\im s_2\subset S^{4n+3}\setminus\{\sqrt{2}p\}$. Hence if $f_1$ is
from $P$ or $A_1\times A_2\times\Sp(n)$ and $f_2\in A_1\times A_2\times\Sp(n)$, then $$F(f_1f_2)=(-e_0)\circ s_1\circ
f_1\circ f_2\circ s_2\circ e_0=(-e_0)\circ s_1\circ f_1\circ s_2\circ e_0\circ (-e_0)\circ s_1\circ f_2\circ s_2\circ
e_0=F(f_1)F(f_2),$$ since $s_2\circ s_1|_{\im s_2}=\id_{\im s_2}$. Therefore it is enough to prove that
$F(f_1f_2)=F(f_1)F(f_2)$, for $f_1\in A_1\times A_2\times\Sp(n)$ and $f_2\in P$. Let
$$f_1=\Op\left(\begin{matrix}a_1e^{-a_2}&0&0\\0&A&0\\0&0&a_1^{-1}e^{-a_2}\end{matrix}\right)\in A_1\times
A_2\times\Sp(n),\,\,
f_2=\Op\left(\begin{matrix}1&-\bar{Y}^{t}&b-\frac{1}{2}Y^{t}\bar{Y}\\0&E_n&Y\\0&0&1\end{matrix}\right)\in P.$$ Then
$f_1f_2f_1^{-1}=f_2'\in P$, where
$f_2'=\Op\left(\begin{matrix}1&-((A^{-1})^t\bar{Y}a_1e^{-a_2})^t&a_1^{2}e^{a_2}(b-\frac{1}{2}Y^{t}\bar{Y})e^{-a_2}\\0&E_n&a_1e^{a_2}(Y^{t}A^{t})^{t}\\0&0&1\end{matrix}\right).$
We have \begin{multline*}F(f_1f_2)=F(f_2'f_1)=F(f_2')F(f_1)=t\left(-\frac{\sqrt{2}}{2}a_1e^{a_2}(Y^{t}A^{t})^{t}\right)
a_1e^{a_2}\Op A\\=t\left(-\frac{\sqrt{2}}{2}a_1e^{a_2}\Op A\cdot Y\right) a_1e^{a_2}\Op A=a_1e^{a_2}\Op A\cdot
t\left(-\frac{\sqrt{2}}{2}Y\right)=F(f_1)F(f_2),\end{multline*} since for any $f\in\mathbb{R}_{+}\times\SO(4n)$ and
$X\in\mathbb{R}^{4n}$ it holds $ft(X)f^{-1}=t(fX)$ or $t(fX)f=ft(X)$. Thus $F$ is the homomorphism from the Lie group
$\Sp(1,n+1)_{\mathbb{H}{p}}$ to the Lie group $\Simil \mathbb{H}^{n}$. Obviously, $F$ is surjective. The claim is proved.

Let $L\subset \mathbb{R}^{4n}$ be a vector (affine) subspace. We call the subset $L\subset \mathbb{H}^{n}$ {\it a real
vector (affine) subspace}.

\begin{theorem}\label{theorem1} Let $G\subset \Sp(1,n+1)_{\mathbb{H}{p}}$ act weakly irreducibly on $\mathbb{H}^{1,n+1}$.
Then if $F(G)\subset \Simil\mathbb{H}^{n}$ preserves a proper real affine subspace $L\subset \mathbb{H}^{n}$, then the
minimal  affine subspace of $\mathbb{H}^{n}$ containing $L$ is $\mathbb{H}^{n}$.
\end{theorem}

{\it Proof.} First we prove that the subgroup $F(G)\subset \Simil \mathbb{H}^{n}$ does not preserve any proper  affine
subspace of $\mathbb{H}^{n}$. Assume that $F(G)$  preserves a  vector subspace $L\subset \mathbb{H}^{n}$. Choosing the
basis $e_1,...,e_n$ of $\mathbb{H}^{n}$ in a proper way, we can suppose that $L=\mathbb{H}^{m}=
\spa_{\mathbb{H}}\{e_1,...,e_m\}$. Consequently, $F(G)\subset
(\mathbb{R}_{+}\times(\Sp(1)\cdot(\Sp(m)\times\Sp(n-m))))\ZR\mathbb{H}^{m}.$ Hence,
$G\subset(\mathbb{R}_{+}\times\Sp(1)\times\Sp(m)\times\Sp(n-m))\ZR(\mathbb{H}^{m}\cdot \real^3)$ and $G$ preserves the
non-degenerate vector subspace $\spa_{\mathbb{H}}\{e_{m+1},...,e_n\}\subset\mathbb{H}^{1,n+1}$. Now suppose that $F(G)$
preserves an affine  subspace $L\subset \mathbb{H}^{n}$. Let $L=Y+L_0$, where $Y\in L$ and $L_0\subset\mathbb{H}^{n}$ is
the vector subspace corresponding to $L$. We may assume that $L_0=\mathbb{H}^{m}=\spa_{\mathbb{H}}\{e_1,...,e_m\}$.
Consider $f=\Op\left(\begin{matrix}1&\sqrt{2}\bar{Y}^{t}&-Y^{t}\bar{Y}\\0&E_n&-\sqrt{2}Y\\0&0&1\end{matrix}\right)\in P$
and the subgroup $\tilde G=f^{-1}Gf\subset\Sp(1,n+1)_{\mathbb{H}{p}}$. For $F(\tilde G)$ we get that $F(\tilde
G)=-t(Y)F(G)t(Y)$. By the above $\tilde G$ preserves the non-degenerate vector subspace
$\spa_{\mathbb{H}}\{e_{m+1},...,e_n\}\subset\mathbb{H}^{1,n+1}$. Hence $G$ preserves the non-degenerate vector subspace
$f(\spa_{\mathbb{H}}\{e_{m+1},...,e_n\})\subset\mathbb{H}^{1,n+1}$. Since $G$ is weakly irreducible, we get $m=n$.

Let $F(G)$ preserve  a real affine subspace $L\subset \mathbb{H}^{n}$ and let $L_0\subset \mathbb{H}^{n}$ be the
corresponding real vector subspace. Consider the vector subspace $(\spa_{\mathbb{H}}L_0)^\bot\subset\mathbb{H}^{n}$. As
above, it can be proved that $G$ preserves the non-degenerate vector subspace
$f((\spa_{\mathbb{H}}L_0)^\bot)\subset\mathbb{H}^{1,n+1}$. Since $G$ is weakly irreducible, we have
$(\spa_{\mathbb{H}}L_0)^\bot=0$ and $\spa_{\mathbb{H}}L_0=\mathbb{H}^n$. The theorem is proved. $\Box$

\section{Proof of the Main Theorem}
First of all, from Example 1 it follows that the algebras of Types I--VIII act weakly irreducibly on
$\mathbb{R}^{4,4n+4}$. For the algebras of Type IX it can be proved in the same way. Therefore we must only prove that any
subalgebra $\g\subset \sp(1,n+1)_{\mathbb{H}{p}}$ that acts weakly irreducibly on $\mathbb{R}^{4,4n+4}$ and contains the
ideal $\B$ is conjugated (by an element from $\SO(4,4n+4))$ to one of the algebras of Types I--IX. Suppose that $\g\subset
\sp(1,n+1)_{\mathbb{H}{p}}$ acts weakly irreducibly on $\mathbb{R}^{4,4n+4}$ and contains the ideal $\B$. Let $G\subset
\Sp(1,n+1)_{\mathbb{H}{p}}$ be the corresponding connected Lie subgroup. By Theorem \ref{theorem1}, $F(G)$ preserves a
real affine subspace $L\subset \mathbb{H}^{n}$ such that the minimal  affine subspace of $\mathbb{H}^{n}$ containing $L$
is $\mathbb{H}^{n}$. We already know that $G$ is conjugated to a subgroup $\tilde G\subset \Sp(1,n+1)_{\mathbb{H}{p}}$
such that $F(\tilde G)$ preserves a real vector subspace $L_0\subset \mathbb{H}^{n}$ with
$\spa_{\mathbb{H}}L_0=\mathbb{H}^n$. Hence we can assume that $F(G)$ preserves a real vector subspace $L\subset
\mathbb{H}^{n}$ and $\spa_\mathbb{H}L=\mathbb{H}^{n}$. Moreover, assume that $F(G)$ does not preserve any proper affine
subspace of $L$. Then $F(G)$ acts transitively on $L$ \cite{Al2}. The connected transitively acting groups of similarity
transformations of the Euclidean spaces are well know. In \cite{GalKahler} these groups were divided into three types. We
describe real subspaces $L\subset \mathbb{H}^{n}$ with $\spa_\mathbb{H}L=\mathbb{H}^{n}$ and subalgebras $\k\subset \lie
(\Simil \mathbb{H}^{n})$ such that the corresponding connected Lie subgroups $K \subset \Simil \mathbb{H}^{n}$ preserve
$L$ and act transitively on $L$. Then the algebra $\g$ must be of the form $(dF)^{-1}(\k)$ for a subalgebra $\k$.

Now we describe real vector subspaces $L\subset \mathbb{H}^{n}$ with $\spa_\mathbb{H}L=\mathbb{H}^{n}.$ Let $L$ be such
subspace. Put $L_1=L\cap iL\cap jL\cap kL$, i.e. $L_1$ is the maximal quaternionic vector subspace in $L$. Let $L_2$ be
the orthogonal complement to $L_1$ in $L$, then $L=L_1\oplus L_2$ and $L_2\cap iL_2\cap jL_2\cap kL_2=0.$ Now let
$L_3=L_2\cap iL_2$, i.e. $L_3$ is the maximal $i$-invariant real vector subspace in $L_2$. Let $L_4$ be its orthogonal
complement in $L_2$, then $L_2=L_3\oplus L_4.$ Similarly, define the spaces $L_5, L_6, L_7, L_8\subset L$ such that
$L_5=L_4\cap jL_4$, $L_4=L_5\oplus L_6$, $L_7=L_6\cap kL_6$ and $L_6=L_7\oplus L_8$. By construction, we get the
orthogonal decomposition  $L=L_1\oplus L_3\oplus L_5\oplus L_7\oplus L_8$ and there exists a $g$-orthogonal basis
$e_1,...,e_n$ of $\mathbb{H}^{n}$ such that this decomposition has the form
\begin{multline}L=\spa_{\mathbb{H}}\{e_1,...e_m\}\oplus \spa_{\mathbb{R}\oplus i\mathbb{R} }\{e_{m+1},...e_{m_1}\}\oplus
\spa_{\mathbb{R}\oplus j\mathbb{R}}\{e_{m_1+1},...e_{m_2}\}\\\oplus \spa_{\mathbb{R}\oplus k\mathbb{R}
}\{e_{m_2+1},...e_{m_3}\}\oplus \spa_{\mathbb{R}}\{e_{m_3+1},...e_n\}.\label{generalL}\end{multline} Obviously, there is
an $f\in \SO(n)$ such that
\begin{equation}\label{fL}fL=\spa_{\mathbb{H}}\{e_1,...e_m\}\oplus \spa_{\mathbb{R}\oplus i\mathbb{R}
}\{e_{m+1},...e_{m+k}\}\oplus \spa_{\mathbb{R}}\{e_{m+k+1},...e_{n}\},\end{equation} where $m+k=m_3$. Since we consider
the subgroups of $\Sp(1,n+1)_{\mathbb{H}{p}}$ up to conjugacy in $\SO(4,4n+4)$, we can assume that $L$ has the form
\eqref{fL}. We will write for short $$L=\mathbb{H}^{m}\oplus \mathbb{C}^{k}\oplus \mathbb{R}^{n-m-k}.$$

Suppose that a subgroup $K\subset \Simil\mathbb{H}^{n}$ preserves $L$. Since $K\subset \Simil\mathbb{H}^{n}\subset
\Simil^0\mathbb{R}^{4n}=(\mathbb{R}_{+}\times \SO(4n))\ZR\mathbb{R}^{4n}$, we have $K\subset (\mathbb{R}_{+}\times
\SO(L)\times \SO(L^{\perp}))\ZR L$. But $K\subset \Simil\mathbb{H}^{n}$, hence $\pr_{\SO(4n)}K\subset \Sp(1)\cdot\Sp(n)$.
Consequently, $\pr_{\SO(4n)}K=\pr_{\Sp(1)\cdot\Sp(n)}K\subset \Sp(1)\cdot\Sp(n)\cap\SO(L)\times\SO(L^{\perp})$. For the
corresponding subalgebra $\k\subset \lie (\Simil \mathbb{H}^{n})$, we have $\pr_{\sp(1)\oplus \sp(n)}\k\subset
\sp(1)\oplus \sp(n)\cap \so(L)\oplus \so(L^{\perp})$. Considering the matrices of the elements of these algebras in the
basis of $\mathbb{R}^{4n}$, we obtain \begin{align*} \sp(1)\oplus \sp(n)\cap
\so(L)\oplus\so(L^{\perp})=\left\{\begin{array}{ll} \sp(1)\oplus \sp(n),& \text{if } m=n;\\ \sp(m)\oplus \u(n-m)\oplus
i\mathbb{R},&\text{if } 0\leq m<n,\\&\quad  n-m-k=0; \\ \sp(m)\oplus \u(k)\oplus \so(n-m-k),&\text{if } 0\leq m<n,\\&
n-m-k\geq 1.
\end{array} \right. \end{align*} The action of the Lie algebras $\u(n-m)$ and $\so(n-m-k)$ on $\mathbb{C}^{n-m}$ and
$\mathbb{R}^{n-m-k}$, respectively, is described in Section \ref{MainTh}.

Let $E$ be a Euclidean space. In \cite{GalKahler} subalgebras $\k\subset \lie (\Simil E)$ corresponding to connected
transitively acting subgroups of $\Simil E$ were divided into the following three types:
\begin{description}
\item{{\bf Type $\mathbb{R}$.}} $\k=(\mathbb{R}\oplus \h)\ltimes E$, where $\h\subset \so(E)$ is a subalgebra.
\item{{\bf Type $\varphi$.}} $\k=\{\varphi(A)+A|A\in\h\}\ltimes E$, where $\h\subset\so(E)$ is a subalgebra,
$\varphi\in\Hom(\h,\mathbb{R})$,
$\varphi|_{\h'}=0$.
\item{{\bf Type $\psi$.}} $\k=\{A+\psi(A)|A\in\h\}\ltimes U$, where we have an orthogonal decomposition $E=W\oplus U$, $\h\subset
\so(W)$ is a subalgebra, $\psi:\h\to W$ is surjective linear map, $\psi|_{\h'}=0$.
\end{description}

Suppose that $m=n$, i.e. $L=\mathbb{H}^{n}$. If $\k$ is of Type $\mathbb{R}$, then $\k=(\mathbb{R}\oplus \h)\ltimes L$,
where $\h\subset \sp(1)\oplus \sp(n)$ is a subalgebra. If $\h\subset \sp(n)$, then $(dF)^{-1}(\k)$ is of  Type II with
$a_2=0$ and $\phi=0$. Let $\h$ have the form $\h_0\oplus \h_1$, where $\h_0\subset \sp(1)$ and $\h_1\subset \sp(n)$. If
$\dim\h_0=1$, then $(dF)^{-1}(\k)$ is of Type II with $\phi=0$ and $\h$ changed to $\h_1$. If $\dim\h_0=2$ or $3$, then
$(dF)^{-1}(\k)$ is of Type I with $\h$ changed to $\h_1$. Suppose that $\h\neq\pr_{\sp(1)}\h\oplus \pr_{\sp(n)}\h$. If
$\h\cap\sp(1)=0$, then $(dF)^{-1}(\k)$ is of Type II with $a_2=0$. Now let $\dim\h\cap\sp(1)=1$ and let $a_2\in
\h\cap\sp(1)$ be a non-zero element. Obviously, $\h=\{A+\phi(A)|A\in\pr_{\sp(n)}\h\}+\mathbb{R}a_2$, where
$\phi:\pr_{\sp(n)}\h\to \sp(1)$ is a homomorphism, $\phi\neq 0$ and $\im\phi\cap \mathbb{R}a_2=0$. For $A+\phi(A)\in\h$,
we have $[A+\phi(A),a_2]=[\phi(A),a_2]\in\h\cap\sp(1)$. Hence, $[\phi(A),a_2]\subset\mathbb{R}a_2$. If $\rk \phi=1$, then
$(dF)^{-1}(\k)$ is of Type II. If $\rk \phi=2$, then there exist $A_1, A_2\in \pr_{\sp(n)}\h$ such that $\phi(A_1),
\phi(A_2)$ and $a_2$ span $\sp(1)$. But this is impossibly, since $\sp(1)'=\sp(1)$. In the same way, if
$\dim\h\cap\sp(1)=2$ and $\h=\{A+\phi(A)\}+(\h\cap \sp(1))$, then $\phi=0$. If $\k=\{\varphi(A)+A|A\in\h\}\ltimes L$ is of
Type $\varphi$, then all $(dF)^{-1}(\k)$ can be obtained from the above, since $\k$ is obtained from $(\mathbb{R}\oplus
\h)\ltimes L$ by twisting between $\h$ and $\mathbb{R}$. We will get that $(dF)^{-1}(\k)$ is of Type III or IV. Let $\k$
be of Type $\psi$, i.e. $\k=\{A+\psi(A)\}\ltimes U$, where $L=W\oplus U$ is an orthogonal decomposition, $\h\subset
\so(W)$ is a subalgebra and $\psi:\h\to W$ is surjective linear map, $\psi|_{\h'}=0$. Since $\h\subset \sp(1)\oplus
\sp(n)$, we have $\h\subset \sp(1)\oplus \sp(n)\cap\so(W)=\sp(W\cap iW\cap jW\cap kW)$. We obtain  Type IX for $m=n$. The
case $m<n$ can be consider similarly. If $\k$ is of Type $\real$, then $\g$ is of Type V or VI. If $\k$ is of Type
$\varphi$, then $\g$ is of Type VII or VIII. If $\k$ is of Type $\psi$, then $\g$ is of Type IX.  The theorem is proved.
$\Box$

\begin{remark}\label{remark1}
It is also possible to classify weakly irreducible subalgebras of $\sp(1,n+1)_{\mathbb{H}{p}}$  containing the ideal $\B$
up to conjugacy by elements of $\Sp(1,n+1)$. For this we should consider in addition the real vector subspace
$L\subset\mathbb{H}^n$ of the form \eqref{generalL} such that at least two of the inequalities $m<m_1<m_2<m_3$ hold. Note
that
\begin{multline*}\sp(1)\oplus \sp(n)\cap \so(L)\oplus\so(L^{\perp})=\sp(\spa_{\mathbb{H}}\{e_1,...e_m\})\oplus
\u(\spa_{\mathbb{R}\oplus i\mathbb{R} }\{e_{m+1},...e_{m_1}\})\\\oplus \u(\spa_{\mathbb{R}\oplus
j\mathbb{R}}\{e_{m_1+1},...e_{m_2}\})\oplus \u(\spa_{\mathbb{R}\oplus k\mathbb{R} }\{e_{m_2+1},...e_{m_3}\})\oplus
\so(\spa_{\mathbb{R}}\{e_{m_3+1},...e_n\}).\end{multline*} We should generalize Type IX assuming that $L$ has the form
\eqref{generalL} and we should in addition add two types of Lie algebras:
\begin{description}
\item{{\bf Type X.}} $\g=\{(a_1,A,X,b)|\, a_1\in\mathbb{R},\,\, A\in\h,\,\, X\in L,\,\, b\in\im\mathbb{H}\},$ where $\h\subset \sp(1)\oplus \sp(n)\cap \so(L)\oplus\so(L^{\perp})$ is a subalgebra.

\item{{\bf Type XI.}} $\g=\{(\varphi(A),A,X,b)|\, A\in\h,\,\, X\in L,\,\, b\in{\im}\mathbb{H}\},$ where $\h\subset \sp(1)\oplus \sp(n)\cap \so(L)\oplus\so(L^{\perp})$ is a subalgebra,
$\varphi\in\Hom(\h,\mathbb{R})$, $\varphi|_{\h'}=0$.
\end{description}
\end{remark}

\bibliographystyle{unsrt}

\vskip1cm

Department of Algebra and Geometry, Masaryk University in Brno, Jan\'a\v ckovo n\'am. 2a, 66295 Brno, Czech Republic

{\it E-mail address}: bezvitnaya@math.muni.cz

\end{document}